\documentclass[12pt,reqno]{amsart}
\usepackage{url}
\usepackage{amsmath,amsthm}
\usepackage{amssymb}
\usepackage[utf8]{inputenc}
\usepackage[T1]{fontenc}
\def\F{\mathbb{F}}

\def\M{\mathcal{M}}

\def\Z{\mathbb{Z}}

\newcommand{\abs}[1]{\left\lvert #1 \right\rvert}

\newtheorem{trm}{Theorem}
\newtheorem{prop}[trm]{Proposition}

\theoremstyle{remark}

\newtheorem{rmk-pl}[rmk]{Remarks}

\newtheorem{ex-sg}[ex]{Example}
\theoremstyle{definition}
\newtheorem*{dfn}{Definition}
\title{Polynomial equations in $\F_q[t]$}
\author{Pierre-Yves Bienvenu}
\date{\today}
\begin{document}
\begin{abstract}
The breakthrough paper of Croot, Lev, Pach \cite{CLP}
on progression-free sets in $\Z_4^n$
introduced a polynomial method that has generated a wealth of
applications, such as Ellenberg and Gijswijt's solutions to the
cap set problem \cite{EG}.
Using this method,
we bound the size of a set of polynomials over $\F_q$ of degree less than $n$ that is free of solutions to the equation
$\sum_{i=1}^k a_if_i^r=0$, where the coefficients $a_i$ 
are polynomials that sum to 0 and the number of variables
satisfies $k\geq 2r^2+1$. The bound 
we obtain
is of the form $q^{cn}$ for some constant $c<1$.
This is in contrast to the best bounds known
for the corresponding problem in the integers, which offer only a logarithmic
saving, but work already with as few as $k\geq r^2+1$ variables.
\end{abstract}
\maketitle
Let $R$ be a ring and $a_1,\ldots,a_k$ be elements of $R$ which sum to 0, i.e. $\sum_{i=1}^k a_i=0$. Then the equation
\begin{equation}\label{eqbase}
\sum_{i=1}^k a_if_i^r=0
\end{equation}
possesses a wealth of trivial solutions $(f_1,\ldots,f_k)$, namely
constant tuples $(f,\ldots,f)$, even though it is not a translation-invariant equation. 
This suggests that if a subset $A\subset R$ is free of non-trivial solutions, then it should be small.
For the ring $R=\Z$, this question was studied first
by Smith \cite{Smith}, Keil \cite{Keil} and Henriot \cite{Henriot}; they replaced the single equation
by a system comprising the initial equation and a linear equation in order to ensure invariance under translation and dilation.  Recently, Browning and Prendiville \cite{BP}
showed that if $r=2$ and $k\geq 5$, and $A\subset [N]$ satisfies $\abs{A}\gg N$ and $N$ is large enough, then 
 equation \eqref{eqbase} necessarily admits
 non-trivial solutions $(f_1,\ldots,f_k)$ in $A^k$. Their method relies on the transference principle. Further, Chow \cite{Chow}  proved that any relatively dense subset of the primes contains a solution to any equation
 of the form \eqref{eqbase}, as long as $k\geq r^2+1$.

Similarly, one may ask whether any dense
subset $A$ of the ring $R=\F_q[t]$ is
bound to contain a non-trivial solution to \eqref{eqbase}. 
In this note,
we answer the question under a natural condition
on the number of variables, namely $k\geq 2r^2+1$.
In the function field setting, the polynomial method of Croot, Lev and Pach \cite{CLP} can be fruitfully applied
and delivers much stronger bounds than any method known in the integers. This was already noticed by Green \cite{Green} in the case of Sarközy's theorem. 

We now precisely state our main theorem. We fix a prime power $q$ and write $P_{q,n}$ for the set of
polynomials of degree strictly less than $n$ over $\F_q$,
so that $\abs{P_{q,n}}=q^n$.

\begin{trm}\label{mytrm}
Let $r,k$ and $d$ be integers
satisfying $k\geq 2r^2+1$. Suppose $(a_1,\ldots,a_k)$ are polynomials 
over $\F_q$ of degree at most $d$ satisfying $\sum_{i=1}^k a_i=0$. Then there exist constants
$0<c(r,q)<1$ and $C=C(d,r,q)$ such that
 any $A\subset P_{q,n}$ satisfying $\abs{A}\geq kC q^{c(r,q)n}$
 must contain a non-trivial solution to the equation \eqref{eqbase}.
 \end{trm}
 
The aforementioned paper of Chow \cite{Chow} implies that $k\geq r^2+1$ is sufficient in the integers, but the bound on the size of $A$ obtained by his analytic method is much weaker 
(we get a power saving, as opposed to his logarithmic saving).
  
We reduce the theorem to the following proposition, which is then tractable by the polynomial method of
Croot-Lev-Pach.
\begin{prop}\label{PolyMap}
For any $\epsilon \in (0,1/2)$, there exists a constant
$c'(\epsilon,q)\in (0,1)$ such that the following holds.
Let $\Phi : (\F_q^n)^k\rightarrow \F_q^m$
be a polynomial map of degree at most $\ell$ (i.e. each coordinate is a polynomial of degree at most $\ell$) and $A\subset \F_q^n$.
Suppose that for any $(f_1,\ldots,f_k)\in A^k$, the equality
$\Phi(f_1,\ldots,f_k)=0$ holds if, and only if, $(f_1,\ldots,f_k)=(f,\ldots, f)$ for some $f\in\F_q^n$.
Finally, suppose that $m\ell/k\leq(1/2-\epsilon)n$. Then $\abs{A}\leq kq^{c'(\epsilon,q)n}$.
\end{prop}

We prove that Proposition \ref{PolyMap} implies Theorem
\ref{mytrm}.
Each polynomial $f=\sum_{i=0}^{n-1}f_it^i\in P_{q,n}$ can be seen as a
vector $\overrightarrow{f}=(f_0,\ldots,f_{n-1})\in\F_q^n$.
Now
$f^r\in P_{q,(n-1)r+1}$
so we see it as the vector $$\overrightarrow{f^r}=(f_0^r,rf_0^{r-1}f_1,\ldots,f_{n-1}^r)\in\F_q^{(n-1)r+1}.$$ We notice that
$\overrightarrow{f^r}=Q(\overrightarrow{f})$ where $Q$ is a polynomial map of degree $r$. Similarly, if $a\in \F_q[t]$
of degree at most $d$, we see that
$f\mapsto \overrightarrow{af}\in\F_q^{n+d}$
is a polynomial map $\F_q^n\rightarrow\F_q^{n+d}$ (of degree 1).
Thus,
$$
\Phi : (f_1,\ldots,f_k)\mapsto \sum_{i=1}^k a_if_i^r
$$
induces a polynomial map of degree $r$
$$
\overrightarrow{\Phi} : (\F_q^n)^k\rightarrow \F_q^{m}
$$
where $m=(n-1)r+d+1$
and
$$
\overrightarrow{\Phi}(\overrightarrow{f_1},\ldots,\overrightarrow{f_k})=\overrightarrow{\Phi(f_1,\ldots,f_k)}.
$$
We observe that if $A\subset P_{q,n}$ does not contain any non-trivial solution to
\eqref{eqbase},
the set 
$\overrightarrow{A}=\{\overrightarrow{f}\mid f\in A\}\subset\F_q^n$ contains only trivial solutions $(\overrightarrow{f},\ldots,\overrightarrow{f})$
to the equation $\overrightarrow{\Phi}(\overrightarrow{f_1},\ldots,\overrightarrow{f_k})=0$.

Moreover, given that 
$k\geq 2r^2+1$, we have
$$\frac{mr}{kn}=\frac{(n-1)r^2+dr+r}{kn}
\leq \frac{r^2}{2r^2+1}+\frac{(d+1)r}{(2r^2+1)n}.$$ 
Hence if
$n\geq 4(d+1)r$, we have $mr/k\leq (1/2-\epsilon)n$,
with $$\epsilon=\epsilon(r)=\frac{1}{4(2r^2+1)}\in (0,1/2).$$
We can then apply Proposition 2 
and obtain $\abs{A}\leq q^{c(r,q)n}$ for some constant
$c(r,q)=c'(\epsilon(r),q)\in (0,1)$. Taking care separately of
the small values of $n$, one can 
find a constant $C(d,r,q)\leq q^{4(d+1)r}$
such that the bound
$$
\abs{A}\leq kC(d,r,q)q^{c(\epsilon,q)n}
$$
is valid for all $n$.

We now prove Proposition \ref{PolyMap}.
We remark, in the spirit of Tao's blog post \cite{Tao}, that the fact that
$$
\forall (f_1,\ldots,f_k)\in A^k,\quad \Phi(f_1,\ldots,f_k)=0\Leftrightarrow f_1=\cdots=f_k
$$
implies that
\begin{equation}
\label{functionEquality}
\forall (f_1,\ldots,f_k)\in A^k,\quad \prod_{i=1}^m
(1-\Phi_i^{q-1}(f_1,\ldots,f_k))=\sum_{f\in A}\prod_{j=1}^k\delta_f(f_j)
\end{equation}
where $\delta_f(f_j)$ is 0 if $f\neq f_j$ and 1 otherwise.
We now recall the notion of slice-rank, as in Tao's blog post or
the article of Kleinberg, Sawin and Speyer
\cite{KSS}.
Take a subset $A\subset \F_q^n$ and
a map $P : A^k\rightarrow \F_q$. Let $\M$ be the set of functions $\F_q^n\rightarrow\F_q$.
This set of functions is naturally in bijection with the set of polynomials in $\F_q[t_1,\ldots,t_n]$
in which no indeterminate is raised to a power greater than $q-1$ (see for instance \cite{Green} for a proof of this bijection).
\begin{dfn}
A \emph{polynomial cover} for $P$ is a tuple $(M_1,\ldots,M_k)\in\M^k$
such that for each $j\in [k]$ and $p\in M_j$, there exists
a function $F_{j,p}$ from  $A^{k-1}$ to  $\F_q$ 
such that for any $(X_1,\ldots,X_k)\in A^k$, we have
\begin{equation}
\label{eqmonomialcover}
P(X_1,\ldots,X_k)=\sum_{j\in[k]}\sum_{p\in M_j}p(X_j)
F_{j,p}(X_1,\ldots,X_{j -1},X_{j+1},\ldots,X_k).
\end{equation}
The \emph{slice rank} of $P$
is the minimum size $\sum_{j\in[k]}\abs{M_j}$ of a polynomial cover.
\end{dfn}

We find that the slice-rank of the right-hand side of \eqref{functionEquality} is $\abs{A}$; indeed, this is
\cite[Lemma 1]{Tao}.
The left-hand side of \eqref{functionEquality} is a polynomial
$P$
in $k\times n$ variables denoted $f_{j,i}$ for
$j\in[k]$ and $i\in[n]$,
and its total degree is at most $(q-1)m\ell$. 
Now $P$ is a sum of monomials of the form
$$
p(f_1,\ldots,f_k)=\prod_{j\in [k]}p_j(f_{j,1},\ldots,f_{j,n})
$$
where each $p_j$ is a monomial in $n$ variables.
For each monomial $p$, by the pigeonhole principle, there exists $j\in[k]$ such that
$\deg p_j\leq (q-1)m\ell/k$.

Being interested in $P$ as a function on $(\F_q^n)^k$,
we reduce it modulo the ideal $I$ generated by the polynomials
$f_{j,i}^{q-1}-f_{j,i}$ for $i\in[n]$ and $j\in[k]$. We continue to use $P$ for the only polynomial in the class $P$ modulo $I$ which 
has degree at most $q-1$ in each variable $f_{j,i}$.
Further, we denote by
$\M_{d,n}$
the set of monomials in $n$ variables, of degree at most $q-1$ in each variable and at most $d$ in total.

We infer from the data above that there exist sets of monomials $M_1,\ldots,M_k\subset \M_{(q-1)m\ell/k,n}$  and
functions $F_{j,p}$ for $(j,p)\in[k]\times M_j$ such that
$$
P=\sum_{j=1}^k\sum_{p\in M_j}p(f_j)F_{j,p}(f_1,\ldots,f_{j-1},f_{j+1},\ldots,f_k).
$$
Now $\abs{\M_{d,n}}/q^n$ may be interpreted as the probability
that the sum of $n$ independent, uniform random variables on $\{0,\ldots,q-1\}$ is at most $d$.
To bound this probability, we use Hoeffding's concentration inequality, which implies that
$$
\abs{\M_{(q-1)m\ell/k,n}}\leq\abs{\M_{(q-1)n(1/2-\epsilon),n}}
\leq q^ne^{-\frac{n\epsilon^2}{2}}=q^{c(\epsilon,q)n}
$$
where
$c(\epsilon,q)=(1-\frac{\epsilon^2}{2\log q})\in(0,1)$.
This implies that the slice-rank of $P$ is at most
$kq^{c(\epsilon,q)n}$ and
concludes the proof of Proposition 2.

In fact, we observe that the scope of our theorem
encompasses more general equations than the diagonal equation \eqref{eqbase}, because Proposition 2 does not require any information on $\Phi$ other than its degree.

\end{document}